\documentclass[12pt]{article}
\usepackage{amsfonts, amsmath}
\usepackage{amsfonts}
\usepackage{amssymb}

\newtheorem{theorem}{Theorem}[section]

\newtheorem{definition}[theorem]{Definition}

\newtheorem{lemma}[theorem]{Lemma}

\newtheorem{proposition}[theorem]{Proposition}
\newtheorem{remark}[theorem]{Remark}

\newenvironment{proof}[1][Proof]{\textbf{#1.} }{\ \rule{0.5em}{0.5em}}
\begin{document}

\author{Tsasa Lusala  \quad \quad J\c{e}drzej \'{S}niatycki  \quad \quad  Jordan Watts}

\title{Regular Points of a Subcartesian Space}

\date{}
\maketitle

\begin{abstract}
We discuss properties of the regular part $S_{reg}$ of a subcartesian space $%
S$. We show that $S_{reg}$ is open and dense in $S$ and the restriction to $%
S_{reg}$ of the tangent bundle space of $S$ is locally trivial.
\end{abstract}


\medskip

\section{Introduction}

In 1967, Sikorski started to study smooth structures on topological spaces in terms
of their corresponding rings of smooth functions \cite{sikorski67}. He
introduced the concept of a \emph{differential space} which is a
generalization of the notion of a smooth manifold. This concept has the advantage that the category of differential spaces is closed under the operation of taking subsets. In
other words, every subset of a differential space inherits a structure of a
differential space such that the inclusion map is smooth, \cite{sikorski71},
\cite{sikorski72}. The theory of differential spaces has been further developed
by several authors.

Also in 1967, Aronszajn introduced the notion of a \emph{subcartesian space}
\cite{aronszajn}, which can be described as a Hausdorff differential space
that is locally diffeomorphic to a differential subspace of a Euclidean
space $\mathbb{R}^{n}$. The original definition of Aronszajn uses a singular
atlas rather than the differential structure provided by the ring of smooth
functions. In the literature on differential spaces, subcartesian spaces as introduced by
Aronszajn are called \emph{differential spaces of class }$D_{0}$, \cite%
{spallek69}, \cite{walczak73}. We use here the term \emph{subcartesian
(differential) spaces} because it is more descriptive.

For a differential space $S,$ with the ring $C^{\infty }(S)$ of smooth
functions on $S$, vectors tangent to $S$ at $x\in S$ are defined as
derivations of $C^{\infty }(S)$ at $x$. They form a vector space denoted by $%
T_{x}S$. The tangent bundle space of $S$ is the set $TS=\cup _{x\in S}T_{x}S$
with an induced structure of a differential space such that the map $\tau
:TS\rightarrow S,$ defined by $\tau ^{-1}(x)=T_{x}S$ for every $x\in S$, is
smooth. It should be mentioned that the dimension of $T_{x}S$ may depend on $x\in S$%
. In the literature, $TS$ is also called the tangent pseudobundle of
$S$, or the Zariski tangent bundle space of $S$ and is denoted by $T^{Z}S$.

A point $x\in S$ is said to be \emph{regular} if there exists a
neighbourhood $U$ of $x$ in $S$ such that $\dim T_{y}S=\dim T_{x}S$ for all $%
y\in U$. Instead of using the dimension of the tangent space $T_{x}S$ at $x$%
, we could use the \emph{structural dimension} of $S$ at $x,$ which is
defined as the minimum of all natural numbers $n$ for which there exists a diffeomorphism of a neighbourhood of $x$ in $S$ onto a subset of $\mathbb{R}^{n}$, \cite%
{marshall75}. We will show that these two notions of dimension are equivalent.

The regular component of a subcartesian space $S$ is the set $S_{reg}$
consisting of regular points of $S$. The aim of this note is to show that,
for a subcartesian space $S$, the regular component $S_{reg}$ is an open and
dense subset of $S$ and that the restriction of $TS$ to $S_{reg}$ is a
locally trivial fibration. It should be noted that $S_{reg}$ need not be a
manifold. For example, for commonly discussed fractals, like the Koch curve
or the Sierpinski gasket, all points are regular. Throughout this paper we
follow the terminology and notations from \cite{sniatycki} and \cite{jordan}.

\section{Preliminaries}

Let $S$ be a subcartesian space, i.e., a Hausdorff
differential space $S$ such that for every point $p\in $ $S,$ there exists $%
n\in \mathbb{N}$ and a neighbourhood of $p$ diffeomorphic to a differential
subspace of $\mathbb{R}^{n}$ which need not be open. For a subcartesian
space $S$, local analysis in a sufficiently small open subset $U$ of $S$ can
be performed in terms of its diffeomorphic image embedded in $\mathbb{R}^{n}$%
. Hence, most of our analysis will be done in terms of differential
subspaces of $\mathbb{R}^{n}$.

Let $S$ be a differential subspace of $\mathbb{R}^{n}$. A function $%
f:S\rightarrow \mathbb{R}$ is smooth if, for every $x\in S$, there exists a
neighbourhood $U$ of $x$ in $\mathbb{R}^{n}$ and a function $f_{x}\in C^{\infty }(%
\mathbb{R}^{n})$ such that
\[
f_{\mid U\cap S}=f_{x\mid U\cap S}.
\]
Thus, the differential structure of $S$ is determined by the ring
\begin{equation}
R(S)=\{f_{|S}\colon f\in C^{\infty }(\mathbb{R}^{n})\}  \label{R(S)}
\end{equation}
consisting of restrictions to $S$ of smooth functions on $\mathbb{R}^{n}$.
Let $N(S)$ denote the ideal of functions in $C^{\infty }(\mathbb{R}^{n})$
which identically vanish on $S$:
\begin{equation}
N(S)=\{f\in C^{\infty }(\mathbb{R}^{n})\colon f_{|S}=0\}.  \label{N(S)}
\end{equation}
We can identify $R(S)$ with the quotient $C^{\infty }(\mathbb{R}^{n})/N(S)$.

\bigskip

Let $S$ be a differential space and $C^{\infty }(S)$ the ring of smooth
functions on $S$. For $x\in S$, a derivation of $C^{\infty }(S)$ at $x$ is a
linear map $u:C^{\infty }(S)\rightarrow \mathbb{R}:f\mapsto u\cdot f$
satisfying Leibniz' rule
\begin{equation}
u\cdot (fh)=(u\cdot f)h(x)+(u\cdot h)f(x)  \label{Leibniz}
\end{equation}%
for every $f,h\in C^{\infty }(S)$. Derivations at $x$ of $C^{\infty }(S)$
form the tangent space\emph{\ }of $S$ at $x$ denoted by $T_{x}S$. The union of
the tangent spaces $T_{x}S$, as $x$ varies over $S$, is the tangent bundle of $S$
and is denoted by $TS$. We denote by $\tau _{S}\colon TS\rightarrow S$ the
tangent bundle projection defined such that $\tau _{S}(u)=x$ if $u\in T_{x}S$%
. The differential structure of the tangent bundle space of a differential
space and smoothness of the tangent bundle projection have been discussed in
\cite{sikorski71}. For the sake of completeness, we describe the
differential structure of $TS$ for a subcartesian space $S$.

Consider first a differential subspace $S$ of $\mathbb{R}^{n}$. We denote by
$q_{1},...,q_{n}$ the restrictions to $S$ of the canonical coordinate
functions\ \thinspace $(x_{1},...,x_{n})$ on $\mathbb{R}^{n}$. For every
function $f\in C^{\infty }(S)$ and $x\in S$, there exists a neighbourhood $U$
of $x$ in $\mathbb{R}^{n}$ and $F\in C^{\infty }(\mathbb{R}^{n})$ such that
\begin{equation}
f_{\mid U\cap S}=F(q_{1},...,q_{n})_{\mid U\cap S}.  \label{1}
\end{equation}
Consider $v\in T_{x}S$, and let $v_{i}=v\cdot q_{i}$ for $i=1,...,n$.
Equation (\ref{1}) yields
\begin{equation}
v\cdot f=(\partial _{1}|_xF)(v\cdot q_{1})+...+(\partial
_{n}|_xF)(v\cdot q_{n})=v_{1}\partial _{1}|_xF+...+v_{n}\partial
_{n}|_xF. \label{2}
\end{equation}
Equation (\ref{2}) shows that $v\in T_{x}S$ can be identified with a vector $%
(v_{1},...,v_{n})\in \mathbb{R}^{n}$. Since $T_{x}S$ has the structure of a
vector space, the set
\begin{equation}
V_x=\{(v_{1},...,v_{n})\in \mathbb{R}^{n}\mid v\mathbf{\in }T_{x}S\}
\label{V}
\end{equation}
is a vector subspace of $\mathbb{R}^{n}$. The tangent bundle $TS$ can be
presented as a subset of $\mathbb{R}^{2n}$ as follows
\begin{equation}
TS=\{(x,v)=(q_{1},...,q_{n},v_{1},...,v_{n})\in \mathbb{R}^{2n}\mid x\in S%
\mbox{
and }v\in V_{x}\}.  \label{TS}
\end{equation}
We denote by $\tau _{S}:TS\rightarrow S$ the tangent bundle projection given
by $\tau _{S}(x,v)=x,$ for every\thinspace $(x,v)\in TS.$

For every $f\in C^{\infty }(S),$ the differential of $f$ is a function $%
df\colon TS\rightarrow \mathbb{R}$ given by
\begin{equation}
df(v)=v\cdot f\mbox{,}  \label{df}
\end{equation}
for every $v\in TS$. The differential structure of $TS$ is generated by the
family of functions $\{q_{1}\circ \tau _{S},...,q_{n}\circ \tau
_{S},dq_{1},...,dq_{n}\}$. In other words, a function $h:TS\rightarrow
\mathbb{R}$ is smooth if, for every $v\in TS$, there is a neighbourhood $W$
of $v$ in $\mathbb{R}^{2n}$ and $H\in C^{\infty }(\mathbb{R}^{2n})$ such
that
\begin{equation}
h_{\mid W\cap TS}=H(q_{1}\circ \tau _{S},...,q_{n}\circ \tau
_{S},dq_{1},...,dq_{n})_{\mid W\cap TS}.  \label{H}
\end{equation}
For $f\in C^{\infty }(S)$ satisfying equation (\ref{1}), we have
\[
f\circ\tau_{S|\tau_S^{-1}(U)}=F(q_1\circ\tau_S,...,q_n\circ\tau_S)_{|\tau_S^
{-1}(U)},
\]
which implies that $\tau _{S}^{\ast }f=f\circ \tau _{S}\in C^{\infty }(TS)$.
Thus, the tangent bundle projection $\tau _{S}$ is smooth.

As before, let $S$ be a differential subspace of $\mathbb{R}^{n}$. A
derivation $v$ of $C^{\infty }(S)$ at $x\in S$ restricts to a derivation of $%
R(S)$ at $x$.

\begin{proposition}
Every derivation of $R(S)$ at $x$ extends to a unique derivation of $%
C^{\infty }(S)$ at $x$. \label{theorem1} \label{propo22}
\end{proposition}

\begin{proof}
Let $w$ be a derivation of $R(S)$ at $x\in S$. Consider $f\in C^{\infty }(S)$%
. There exist an open neighbourhood $U$ of $x$ in $\mathbb{R}^{n}$ and a
function $f_{x}\in C^{\infty }(\mathbb{R}^{n})$ such that $f_{|U\cap S}={%
f_{x}}_{|U\cap S}$. Set $\widetilde{w}(f)=w({f_{x}}_{|S})$. Let $V$
be another open neighbourhood of $x$ in $\mathbb{R}^{n}$ and
$g_{x}\in C^{\infty }(\mathbb{R}^{n})$ a function such that
$f_{|V\cap S}={g_{x}}_{|V\cap S}$. We have that $U\cap V\cap S$ is
an open subset of $S$ and ${f_{x}}_{|U\cap V\cap S}={g_{x}}_{|U\cap
V\cap S}$. Therefore $({f_{x}-g_{x}})_{|U\cap V\cap S}=0$, i.e.,
$({f_{x}-g_{x}})_{|S}\in R(S)\subset C^{\infty }(\mathbb{R}^{n})
$ vanishes identically on the open subset $U\cap V\cap S$ of $S$. Hence, $w({%
f_{x}}_{|S}-{g_{x}}_{|S})=0$. This proves that the extension $\widetilde{w}$
is a well-defined derivation of $C^{\infty }(S)$ extending the derivation $w$
of $R(S)$ at $x$. Finally, it is clear that such an extension $\widetilde{w}$
of $w$ is uniquely defined.
\end{proof}

\begin{remark}
Equation (\ref{2}) shows that every derivation of $C^{\infty }(S)$ at $x\in
S\subseteq \mathbb{R}^{n}$ can be extended to a derivation of $C^{\infty }(%
\mathbb{R}^{n}).$ We can ask the question under what conditions a
derivation $w$ of $C^{\infty }(\mathbb{R}^{n})$ at $x\in S\subseteq
\mathbb{R}^{n}$ defines a derivation of $C^{\infty }(S)$ at
$x$.\label{remark1}
\end{remark}

\begin{proposition}
A derivation $w$ of $C^{\infty }(\mathbb{R}^{n})$ at $x\in S\subseteq
\mathbb{R}^{n}$ defines a derivation of $C^{\infty }(S)$ at $x$ if and only
if $w$ annihilates $N(S)$, i.e. $w(f)=0$ for all $f\in N(S)$. \label%
{propo21}
\end{proposition}

\begin{proof}
It follows from Proposition \ref{theorem1} and Remark \ref{remark1} that
derivations at $x$ of $C^{\infty }(S)$ can be identified with derivations at
$x$ of $R(S)$. Now one uses the identification
\[
R(S)\equiv \frac{C^{\infty }(\mathbb{R}^{n})}{N(S)}=\frac{C^{\infty }(%
\mathbb{R}^{n})}{\sim {}},
\]%
where $f\,{\sim }\,g$ in $C^{\infty }(\mathbb{R}^{n})$ if and only if $%
f-g\in N(S)$. For a derivation $w$ at $x$ of $C^{\infty }(\mathbb{R}^{n})$,
one defines $w([f])=w(f)$. It is clear that this defines a derivation of $%
R(S)$ if and only if $w(f)=0$ for all $f\in N(S)$.
\end{proof}

\bigskip

\section{The regular component of a subcartesian space}

We now discuss the notion of structural dimension introduced by
Marshall, \cite{marshall75}.

\begin{definition}
Let $S$ be a subcartesian space. The structural dimension of a point
$x\in S$ is the smallest integer, denoted by $n_x$, such that for
some open
neighbourhood $U\subseteq S$ of $x$, there is a diffeomorphism $%
\varphi:U\rightarrow V$ for some arbitrary subset $V\subseteq \mathbb{R}^{n}$%
.
\end{definition}

A real-valued function $f:D\rightarrow \mathbb{R}$ is upper
semi-continuous if the subset of $D$ determined by $\{x \in
D\colon{} f(x)< a\}$, for any $a\in\mathbb{R}$, is open.

\begin{lemma}
The function $N\colon S\rightarrow\mathbb{N}:x\mapsto n_x $ is upper
semi-continuous.
\end{lemma}

\begin{proof}
Let $S_i = \{x\in S\colon{} n_x \leq i\}$. Assume that $S_i$ is not open.
Then there exists a point $z\in S_i$ such that there is no open
neighbourhood $U\subseteq S_i$ of $z$. But then, there is no open
neighbourhood $V\subseteq S$ of $z$ diffeomorphic to an arbitrary subset of $%
\mathbb{R}^j$ for any $j\leq i$. Hence, $n_z>i$, and so $z$ is not in $S_i$.
Thus, $S_i$ is open, and so the structural dimension serves as an upper
semi-continuous function on $S$.
\end{proof}

\begin{definition}
A point $x\in S$ is called a structurally regular point if there is
a neighbourhood $U$ of $x$ in $S$ such that $n_y=n_x $ for all $y\in
U$.  A point that is not structurally regular is called structurally
singular.
\end{definition}

The regular component $S_{reg}$ of a subcartesian space $S$ is the
set of all structurally regular points of $S.$

\begin{lemma}
For every point $x$ of a subcartesian space $S$, the structural dimension of
$S$ at $x$ is equal to $\dim T_{x}S.$
\end{lemma}

\begin{proof}
Let $n=n_x$. So, there is a neighbourhood $U\subseteq S$ of $x$
diffeomorphic to a differential subspace of $\mathbb{R}^n$. Since any
derivation of $C^{\infty}(S)$ can be extended to a derivation of $C^{\infty}(%
\mathbb{R}^n)$, we have $\mbox{dim}\,T_xS\leq\mbox{dim}\,\mathbb{R}^n=n$.

Now assume that $\mbox{dim}\,T_{x}S<n$, then there exists a
derivation $u\in T_{x}\mathbb{R}^{n}$ that is not an extension of a
derivation of $C^{\infty
}(S)$. This implies by Proposition \ref{propo21} that there is a function $%
f\in N(U)$ such that $u(f)\not=0$. In this case, if $p^{1}$, ..., $p^{n}$
are the canonical coordinate functions on $\mathbb{R}^{n}$, then ${\partial
_{p^{j}}}f_{\mid x}\not=0$, for some $j\in \{1,...,n\}$. Hence, there is a
neighbourhood $V\subseteq f^{-1}(0)$ of $x$ that is a submanifold of $%
\mathbb{R}^{n}$. It is clear that the structural dimension of $S$ at
points in $V$ is $m<n$ ($m$ being the dimension of $V$ as a
manifold). There exists an open
neighbourhood $\tilde{V}\in V$ of $x$ diffeomorphic to an open subset of $%
\mathbb{R}^{m}$. Since $f\in N(U)$, there exists a neighbourhood $W\subset
U\subset f^{-1}(0)$ of $x$. So $\tilde{V}\cap W$ is a neighbourhood of $x$
in $\mathbb{R}^{m}$. This is a contradiction as the structural dimension $%
n_{x}=n>m$. Therefore, $\dim \,T_{x}S=n_{x}$.
\end{proof}

\begin{lemma}
Let $n$ be the maximum of the structural dimensions of $S$ at points
of an open subset $V\subset S$. If every open subset contained in
$V$ has a point at which the structural dimension is $n$, then $V$
consists of regular points. \label{denselemma}
\end{lemma}

\begin{proof}
The assumption implies that the subset $W=\{x\in V\colon n_{x}=n\}$
is dense in $V$. For each $x\in V$, let $O_{x}$ be an open
neighbourhood of $x$ in $V$ diffeomorphic to a subset of
$\mathbb{R}^{n}$. Take $y\in V\setminus W$. Then $n_{y}<n$ (by the
definition of $n$). Let $O_{y}$ be an open neighbourhood of $y$ in
$V$ diffeomorphic to a subset of $\mathbb{R}^{n_{y}}$. Since $W$ is
dense in $V$, there exists $x\in W\cap O_{y}$. So, $O_{x}\cap O_{y}$
is diffeomorphic to a subset of $\mathbb{R}^{n_{y}}$. But $n$ is the
minimum of all $m$ such that a neighbourhood of $x$ is diffeomorphic
to a subset of $\mathbb{R}^{m}$. Since $O_{x}\cap O_{y}$ is a
neighbourhood of $x$ diffeomorphic to a subset of
$\mathbb{R}^{n_{y}}$, we have $n\leq n_{y}$. But $n_{y}<n$ by
assumption. Therefore, $V\setminus W$ is empty, i.e., the dimension
of $S$ at a point of the open subset $V$ is $n$. This implies that
every point in $V$ is structurally regular.
\end{proof}

\begin{theorem}
The set $S_{reg}$ of all structurally regular points of a subcartesian space
$S$ is open and dense in $S$.
\end{theorem}

\begin{proof}
Let $x\in S_{reg}$. Since $x$ is a structurally regular point, there exists
an open neighbourhood $U\subseteq S$ of $x$ such that for every $y\in U$, $%
n_y=n_x$. This implies that every point of $U$ is structurally regular.
Hence, $U\subseteq S_{reg}$. Therefore, $S_{reg}$ is an open subset of $S$.

\noindent Now, suppose that the subset $S_{reg}$ of structural
regular points is not dense in $S$. In this case, there exists a
non-empty open subset $U\subseteq S$ such that $U$ contains no
structurally regular points, i.e., every point in $U$ is a
structurally singular point. Without loss of generality, we assume
that $U$ is diffeomorphic to a differential subspace of
$\mathbb{R}^{n}$ for some $n>0$. In fact, $n$ cannot be $0$,
otherwise $U$ would be a set of isolated points which are regular by
the induced topology. Define $S_{i}=\{x\in S\colon n_{x}\leq i\}$.
Assume that $U\subset S_{k}$ (for some $k>0$). It follows that, if
$V_{1}\subset U$ is an open subset, then $V_{1}$ contains infinitely
many points where the structural dimensions are at least two
different numbers from $0$ to $k$. Let $n_{1}$ be the
maximum of these structural dimensions at points in $V_{1}$. By Lemma \ref%
{denselemma}, there exists an open subset $V_{2}\subset V_{1}$ such that the
maximum of structural dimensions of $S$ at points in $V_{2}$ is $n_{2}<n_{1}$%
. Similarly, there exists an open subset $V_{3}\subset V_{2}$ with a maximum
of structural dimensions at its points, $n_{3}<n_{2}$. Thus, continuing this
process, we have the following decreasing sequence
\[
n_{1}>n_{2}>n_{3}>\cdots> n_{i},
\]%
stopping at some $n_{i}\geq 0$. We reach some open subset
$V_{i}\subset U$ such that the structural dimension at all points of
$V_{i}$ is $n_{i}\geq 0$. Hence, all points of $V_{i}$ are regular
points. As a consequence, since $U$ contains no regular points, $U$
is not a subspace of $S_{k}$ for any $k\geq 0 $. But we are dealing
only with finite structural dimensions and $U$ was chosen to be
diffeomorphic to a differential subspace of $\mathbb{R}^{n}$ for
some $n$, so we have $U\subset S_{n}$, which is a contradiction.
Therefore, a non-empty open subset $U\subset S$ containing no
structurally regular points does not exist. This completes the proof
that the set $S_{reg}$
of all structurally regular points of a subcartesian space $S$ is dense in $%
S $.
\end{proof}

\begin{theorem}
Let $S$ be a subcartesian space. Then the restriction of the tangent bundle
projection $\tau\colon TS\longrightarrow S$ to $T(S_{reg})$ is a locally
trivial fibration over $S_{reg}$. For each $x\in S_{reg}$ with structural
dimension $n$, there is a neighbourhood $W$ of $x$ in $S$ and a family $X_1$%
,...,$X_n$ of global derivations of $C^{\infty}(S)$ such that $%
T_WS=\tau^{-1}(W)$ is spanned by the restrictions $X_1$,...,$X_n$ to $V$.
\end{theorem}

\begin{proof}
Let $x\in S_{reg}$ with $n_x=n$. Since $S_{reg}$ is open, there
exists
a neighbourhood $V\subset S_{reg}$ of $x$ such that $n_y=n$ for all $y \in V$%
. As $S$ is a subcartesian space, we may assume without loss of
generality that there is an embedding $\varphi$ of $V$ into
$\mathbb{R}^n$. We first prove that $TV$ the set of all pointwise
derivations of $C^{\infty}(V)$ is a trivial bundle.

\noindent

Let $R(V)$ consist of restrictions to $V$ of all smooth functions on $%
\mathbb{R}^{n}$, and $N(V)$ be the space of functions on
$\mathbb{R}^{n}$ which vanish on $V$. We identify $R(V)$ with
$C^{\infty }(\mathbb{R}^{n})$ modulo $N(V)$. It follows that
$\partial_{i}|_yf_{|V}=0$ for every $i=1,...,n$, each $f\in N(V)$
and $y\in V\subset\mathbb{R}^n$. By Proposition \ref{propo21}, we
have that $\partial _{1}|_y,...,\partial _{n}|_y$ define derivations
of $C^{\infty }(V)$ at each $y\in V$. Hence, there are $n$ sections
$X_{1}$,..., $X_{n}$ of the tangent bundle projection $\tau
_{V}\colon TV\longrightarrow V$ such that
\[
X_{i}|_y(h\,\,\mbox{mod}\,\,N(V))=(\partial _{i}|_yh)
\]%
for every $i=1,...,n$, $h\in R(V)$ and $y\in V$. Now we need to
prove that the sections $X_{1}$,..., $X_{n}$ are smooth. Let
$q_{1},...,q_{n}$ be
restrictions to $V$ of the coordinate functions on $\mathbb{R}^{n}$. For $%
i=1,...,n$, we denote by $dq_{i}$ the function on $TV$ such that
\[
dq_{i}(w)=w(q_{i})
\]%
for every $w\in TV$. The differential structure of $TV$ is generated by the
functions $(\tau _{V}^{\ast }q_{1},...,\tau _{V}^{\ast
}q_{n},dq_{1},...,dq_{n})$ in the sense that every function $f\in C^{\infty
}(TV)$ is of the form
\[
f=F(\tau _{V}^{\ast }q_{1},...,\tau _{V}^{\ast }q_{n},dq_{1},...,dq_{n})
\]%
for some $F\in C^{\infty }(\mathbb{R}^{2n})$. In order to show that $%
X_{i}\colon V\rightarrow TV$ is smooth, it suffices to show that, for every $%
f\in C^{\infty }(TV)$, the pull-back $X_{i}^{\ast }f$ is in $C^{\infty }(V)$%
. Since
\[
dq_{i}\circ X_{j}=\delta _{ij}=\left\{
\begin{array}{l}
1\quad \mbox{if}\quad i=j \\
0\quad \mbox{if}\quad i\not=j%
\end{array}%
,\right.
\]%
it follows that
\begin{eqnarray*}
X_{i}^{\ast }f &=&f\circ X_{i}=F(\tau _{V}^{\ast }q_{1},...,\tau _{V}^{\ast
}q_{n},dq_{1},...,dq_{n})\circ X_{i} \\
&=&F(\tau _{V}^{\ast }q_{1}\circ X_{i},...,\tau _{V}^{\ast }q_{n}\circ
X_{i},dq_{1}\circ X_{i},...,dq_{n}\circ X_{i}) \\
&=&F(q_{1}\circ \tau _{V}\circ X_{i},...,q_{n}\circ \tau _{V}\circ
X_{i},\delta _{1i},...,\delta _{ni}) \\
&=&F(q_{1},...,q_{n},\delta _{1i},...,\delta _{ni}).
\end{eqnarray*}%
Hence $X_{i}^{\ast }f$ is in $C^{\infty }(V)$. This implies that the tangent
bundle space $TV$ is globally spanned by $n$ linearly independent smooth
sections $X_{1},...,X_{n}$. Thus, $TV$ is a trivial bundle. We can choose an
open neighbourhood $W$ of $x$ contained in $V$ such that its closure $\bar{W}
$ is also in $V$. Using bump functions that are equal to $1$ on $W$ and $0$ outside of $V$, we can construct derivations of $%
C^{\infty}(S)$ that extend restrictions of $X_1,...,X_n$ to $W$.
Hence $TW$
is spanned by the restrictions to $W$ of global derivations of $%
C^{\infty}(S) $. This completes the proof.
\end{proof}

\begin{tabular}{llll}
Tsasa Lusala and J\c{e}drzej \'{S}niatycki  &  & & Jordan Watts \\
Department of Mathematics and  Statistics &  & &Department of Mathematics \\
University of Calgary &  & & University of Toronto \\
2500 University Drive NW &  & &Toronto, Ontario, Canada \\
Calgary, Alberta, Canada &  &  &M5S 2E4 \\
T2N 1N4 &  & &email: jwatts@math.toronto.edu \\
email:\{tsasa, sniat\}@math.ucalgary.ca &  &  &\\
\end{tabular}

\end{document}